\documentclass[11pt]{article}
\usepackage{amsfonts,amsmath,amssymb,latexsym}
\usepackage{eucal}
\usepackage{geometry}

\newtheorem{proposition}{Proposition}
\newtheorem{theorem}{Theorem}
\newtheorem{lemma}{Lemma}
\newtheorem{remark}{Remark}

\begin{document}

\title{Existence of nonparametric solutions for a capillary problem in warped products}
\author{Jorge H. Lira, Gabriela A. Wanderley}
\maketitle

\begin{abstract}
We prove that there exist solutions for a non-parametric capillary problem in a wide class of Riemannian manifolds endowed with a Killing vector field. In other terms, we prove the existence of Killing graphs with prescribed mean curvature and prescribed contact angle along its boundary. These results may be useful for modelling stationary hypersurfaces under the influence of a non-homogeneous gravitational field defined over an arbitrary Riemannian manifold. 

\vspace{3mm}

\vspace{2mm}

\end{abstract}

\section{Introduction}

Let $M$ be a $(n+1)$-dimensional Riemannian manifold endowed with a Killing vector field $Y$. Suppose that the distribution orthogonal to $Y$ is of constant rank and integrable. Given an integral leaf $P$ of that distribution, let $\Omega\subset P$ be a bounded domain with regular boundary $\Gamma =\partial\Omega$. We suppose for simplicity that $Y$ is complete. In this case, let $\vartheta: \mathbb{R}\times \bar\Omega \to M$ be the flow generated by $Y$ with initial values in $M$. In geometric terms, the ambient manifold is a warped product $M = P\times_{1/\sqrt{\gamma}} \mathbb{R}$ where $\gamma = 1/|Y|^2$.

The Killing graph of a differentiable function  $u:\bar\Omega\to \mathbb{R}$  is the hypersurface $\Sigma \subset M$ parametrized by the map
\[
X(x)=\vartheta(u(x),x), \quad x\in\bar\Omega.
\]
The Killing cylinder $K$ over $\Gamma$ is by its turn defined by
\begin{equation}
K=\{\vartheta(s,x): s\in \mathbb{R}, \, x \in \Gamma\}.
\end{equation}
The height function with respect to the leaf $P$ is measured by the arc lenght parameter $\varsigma$ of the flow lines of $Y$, that is, 
\[
\varsigma=\frac{1}{\sqrt\gamma}s. 
\]
Fixed these notations, we are able to formulate a capillary problem in this geometric context which model stationary graphs under a gravity force whose intensity depends on the point in the space. More precisely, given a \emph{gravitational potential} $\Psi \in C^{1,\alpha}(\bar\Omega \times \mathbb{R})$ we define the functional
\begin{equation}
\mathcal{A}[u] = \int_\Sigma \bigg(1+\int_0^{u/\sqrt\gamma}\Psi(x, s(\varsigma)) \,\textrm{d}\varsigma\bigg)\textrm{d}\Sigma. 
\end{equation}
The volume element $\textrm{d}\Sigma$  of $\Sigma$ is given by
\[
\frac{1}{\sqrt\gamma}\sqrt{\gamma+|\nabla u|^2}\,\textrm{d}\sigma,
\]
where $\textrm{d}\sigma$ is the volume element in $P$. 

The first variation formula of this functional may be deduced as follows. Given an aarbitrary function $v\in C^\infty_c(\Omega)$ we compute
\begin{eqnarray*}
& & \frac{d}{d\tau}\Big|_{\tau=0}\mathcal{A}[u+\tau v] =\int_\Omega \bigg(\frac{1}{\sqrt\gamma}\frac{\langle \nabla u, \nabla v\rangle}{\sqrt{\gamma+|\nabla u^2|}} + \frac{1}{\sqrt\gamma}\Psi (x, u(x)) v\bigg) \sqrt{\sigma}\textrm{d}x\\
& & \,\, = \int_\Omega \bigg(\textrm{div}\Big(\frac{1}{\sqrt\gamma}\frac{\nabla u}{W}v\Big) - \textrm{div}\Big(\frac{1}{\sqrt\gamma}\frac{\nabla u}{W}\Big) v + \frac{1}{\sqrt\gamma}\Psi (x, u(x)) v\bigg) \sqrt{\sigma}\textrm{d}x
\\
& & \,\,\,\, -\int_\Omega \bigg(\frac{1}{\sqrt\gamma}\textrm{div}\Big(\frac{\nabla u}{W}\Big) - \frac{1}{\sqrt\gamma}\langle \frac{\nabla \gamma}{2\gamma}, \frac{\nabla u}{W}\rangle  -\frac{1}{\sqrt\gamma}\Psi (x, u(x)) \bigg) v \sqrt{\sigma}\textrm{d}x,
\end{eqnarray*} 
where $\sqrt\sigma \textrm{d}x$ is the volume element $\textrm{d}\sigma$ expressed in terms of local coordinates in $P$. The differential operators $\textrm{div}$ and $\nabla$  are respectively the  divergence and gradient in $P$ with respect to the metric induced from $M$.

We conclude that stationary functions satisfy the capillary-type equation
\begin{equation}
\label{capillary}
\textrm{div}\Big(\frac{\nabla u}{W}\Big) - \langle \frac{\nabla \gamma}{2\gamma}, \frac{\nabla u}{W}\rangle = \Psi. 
\end{equation}
Notice that a Neumann boundary condition arises naturally from this variational setting: given a $C^{2,\alpha}$ function $\Phi:K \to (-1,1)$, we impose the following prescribed angle condition
\begin{equation}
\label{neumann-condition}
\langle N, \nu\rangle = \Phi
\end{equation}
along $\partial\Sigma$, where 
\begin{equation}
N = \frac{1}{W}\big(\gamma Y - \vartheta_* \nabla u\big)
\end{equation}
is the unit normal vector field along $\Sigma$ satisfying $\langle N, Y\rangle >0$ and $\nu$ is the unit normal vector field along $K$ pointing inwards the Killing cylinder over $\Omega$. 

Equation (\ref{capillary}) is the prescribed mean curvature equation for Killing graphs.  A general existence result for solutions of the Dirichlet problem for this equation may be found in \cite{DHL}. There the authors used local perturbations of the Killing cylinders as barriers for obtaining height and gradient estimates. However this kind of barrier is not suitable to obtain \emph{a priori} estimates for solutions of Neumann problems. For that reason we consider now local perturbations of the graph itself adapted from the original Korevaar's approach in \cite{korevaar} and its extension by M. Calle e L. Shahriyari \cite{calle}. 

Following \cite{calle} and \cite{korevaar} we suppose that the data $\Psi$ and $\Phi$ satisfy   
\begin{itemize}
\item[i.] $|\Psi|+|\bar\nabla\Psi|\le C_\Psi$ in $\bar\Omega\times \mathbb{R}$, 
\item[ii.] $\langle \bar\nabla \Psi, Y\rangle \ge \beta>0$ in $\bar\Omega\times \mathbb{R}$,
\item[iii.] $\langle \bar\nabla\Phi, Y\rangle \le 0$,
\item[iv.] $(1-\Phi^2)\ge \beta'$,
\item[v.] $|\Phi|_2\le C_\Phi$ in $K$,
\end{itemize}
for some positive constants $C_\Psi, C_\Phi, \beta$ and $\beta'$, where $\bar\nabla$ denotes the Riemannian connection in $M$. Assumption ($ii)$ is classically referred to as the \emph{positive gravity} condition. Even in the Euclidean space, it seems to be an essential assumption in order to obtain \emph{a priori} height estimates. A very geometric discussion about this issue may be found at \cite{concus-finn}. Condition ($iii$) is the same as in \cite{calle} and \cite{korevaar} since at those references $N$ is chosen in such a way that $\langle N, Y\rangle >0$.  
 
The main result in this paper is the following one

\begin{theorem}
\label{main} Let $\Omega$ be a bounded $C^{3,\alpha}$ domain in $P$. 
Suppose that the $\Psi\in C^{1,\alpha}(\bar\Omega\times\mathbb{R})$ and $\Phi\in C^{2,\alpha}(K)$ with $|\Phi|\le 1$ satisfy conditions {\rm (i)-(v)} above. Then there exists a unique solution $u\in C^{3,\alpha}(\bar\Omega)$ of the capillary problem {\rm (\ref{capillary})-(\ref{neumann-condition})}. 
\end{theorem}

We observe that $\Psi=nH$, where $H$ is the mean curvature of $\Sigma$ calculated with respect to $N$. Therefore Theorem \ref{main} establishes the existence of Killing graphs with prescribed mean curvature $\Psi$ and prescribed contact angle with $K$ along the boundary. Since the Riemannian product $P\times \mathbb{R}$ corresponds to the particular case where $\gamma=1$, our result extends the main existence theorem in \cite{calle}. Space forms constitute other important examples of the kind of warped products we are considering. In particular,  we encompass the case of Killing graphs over totally geodesic hypersurfaces in the hyperbolic space $\mathbb{H}^{n+1}$.

In Section \ref{section-height}, we prove \emph{a priori} height estimates for solutions of (\ref{capillary})-(\ref{neumann-condition}) based on  Uraltseva's method as presented in \cite{uraltseva}. These height estimates are one of the main steps for using the well-known Continuity Method in order to prove Theorem \ref{main}. At this respect, we refer the reader to the classical references \cite{concus-finn}, \cite{gerhardt} and \cite{spruck-simon}.

Section \ref{section-gradient} contains the proof of  interior and boundary gradient estimates. There we follow closely a method due to N. Korevaar \cite{korevaar} for graphs in the Euclidean spaces and extended by M. Calle and L. Shahriyari \cite{calle} for Riemannian products. Finally the classical Continuity Method is applied to (\ref{capillary})-(\ref{neumann-condition}) in Section \ref{section-proof} for proving the existence result. 



%

\section{Height estimates}
\label{section-height}

In this section, we use a technique developed by N. Uraltseva \cite{uraltseva} (see also \cite{uraltseva-book} and \cite{GT} for classical references on the subject) in order to obtain a height estimate for solutions of the capillary problem (\ref{capillary})-(\ref{neumann-condition}). This estimate requires the \emph{positive gravity} assumption ($ii$) stated in the Introduction. 

\begin{proposition} Denote 
\begin{equation}
\beta = \inf_{\Omega\times \mathbb{R}}\langle \bar\nabla \Psi, Y\rangle
\end{equation}
and 
\begin{equation}
\mu = \sup_\Omega \Psi(x,0). 
\end{equation}
Suppose that $\beta >0$. Then any solution $u$  of (\ref{capillary})-(\ref{neumann-condition}) satisfies
\begin{equation}
|u(x)|\le \frac{\sup_\Omega |Y|}{\inf_\Omega |Y|}\frac{\mu}{\beta}
\end{equation}
for all $x\in \bar\Omega$.
\end{proposition}

\noindent \emph{Proof.}
Fix an arbitrary real number $k$ with
\begin{equation*}
k > \frac{\sup_\Omega |Y|}{\inf_\Omega |Y|}\frac{\mu}{\beta}.
\end{equation*}
Suppose that the superlevel set
\begin{equation*}
\Omega_k = \{x\in \Omega: u(x)>k\}
\end{equation*}
has a nonzero Lebesgue measure.  Define $u_k:\Omega \to \mathbb{R}$  as
\begin{equation*}
u_k(x) = \max\{u(x)-k,0\}.
\end{equation*}
From the variational formulation we have
\begin{eqnarray*}
0 &=&\int_{\Omega_k} \bigg(\frac{1}{\sqrt\gamma}\frac{\langle \nabla u, \nabla u_k\rangle}{\sqrt{\gamma+|\nabla u^2|}} + \frac{1}{\sqrt\gamma}\Psi (x, u(x)) u_k\bigg) \sqrt{\sigma}\textrm{d}x\\
&=&  \int_{\Omega_k} \bigg(\frac{1}{\sqrt\gamma}\frac{|\nabla u|^2}{W} +\frac{1}{\sqrt\gamma} \Psi (x, u(x)) (u-k)\bigg) \sqrt{\sigma}\textrm{d}x\\
& = &  \int_{\Omega_k} \bigg(\frac{1}{\sqrt\gamma}\frac{W^2-\gamma}{W} +\frac{1}{\sqrt\gamma} \Psi (x, u(x)) (u-k)\bigg) \sqrt{\sigma}\textrm{d}x \\
&= &
\int_{\Omega_k} \bigg(\frac{W}{\sqrt\gamma}-\frac{\sqrt\gamma}{W} + \frac{1}{\sqrt\gamma}\Psi (x, u(x)) (u-k)\bigg) \sqrt{\sigma}\textrm{d}x  .
\end{eqnarray*}
However
\begin{equation*}
\Psi(x,u(x)) = \Psi(x,0) +\int_0^{u(x)} \frac{\partial \Psi}{\partial s}\textrm{d}s \ge -\mu +\beta u(x). 
\end{equation*}
Since $\frac{\sqrt{\gamma}}{W}\leq 1$ we conclude that
\begin{eqnarray*}
|\Omega_k|-|\Omega_k|-\mu\int_{\Omega_k}\frac{1}{\sqrt{\gamma}}(u-k)+\beta\int_{\Omega_k}\frac{1}{\sqrt{\gamma}}u(u-k)\le 0.
\end{eqnarray*}
Hence  we have
\begin{eqnarray*}
\beta\int_{\Omega_k}\frac{1}{\sqrt{\gamma}}u(u-k) \le \mu\int_{\Omega_k}\frac{1}{\sqrt{\gamma}}(u-k).\nonumber
\end{eqnarray*}
It follows that
\begin{eqnarray*}
\beta k \inf_\Omega |Y| \int_{\Omega_k}(u-k) \le \mu\sup_\Omega |Y|\int_{\Omega_k}(u-k)\nonumber
\end{eqnarray*}
Since $|\Omega_k|\neq 0$ we have 
\[
k \le \frac{\sup_\Omega |Y|}{\inf_\Omega |Y|}\frac{\mu}{\beta},
\]
what contradicts the choice of $k$. We conclude that $|\Omega_k|=0$ for all $k \ge \frac{\sup_\Omega |Y|}{\inf_\Omega |Y|}\frac{\mu}{\beta}$. This implies that 
\[
u(x)\le \frac{\sup_\Omega |Y|}{\inf_\Omega |Y|}\frac{\mu}{\beta},
\]
for all $x\in \bar\Omega$. A lower estimate may be deduced in a similar way.  This finishes the proof of the Proposition. $\hfill\square$

\begin{remark}
The construction of geometric barriers similar to those ones in \cite{concus-finn} is also possible at least in the case where $P$ is endowed with a rotationally invariant metric and $\Omega$ is contained in a normal neighborhood of a pole of $P$. 
\end{remark}

\section{Gradient estimates}
\label{section-gradient}

Let $\Omega'$ be a subset of $\Omega$ and define
\begin{equation}
\Sigma'= \{\vartheta(u(x),x): x\in \Omega'\}\subset \Sigma
\end{equation}
be the graph of $u|_{\Omega'}$. Let $\mathcal{O}$ be an open subset in $M$ containing $\Sigma'$. We consider a vector field $Z\in \Gamma(TM)$ with bounded $C^2$ norm and supported in $\mathcal{O}$. Hence there exists $\varepsilon>0$ such that the local flow  
$\Xi:(-\varepsilon, \varepsilon)\times \mathcal{O}\to M$  generated by $Z$ is well-defined. We also suppose that 
\begin{equation}
\label{Zboundary}
\langle Z(y), \nu (y)\rangle = 0,
\end{equation}
for any $y\in K\cap\mathcal{O}$. This implies that the flow line of $Z$ passing through a point $y\in K\cap\mathcal{O}$ is entirely contained in $K$. 

We define a variation of $\Sigma$ by a one-parameter family of hypersurfaces $\Sigma_\tau$, $\tau \in (-\varepsilon, \varepsilon)$, parameterized by $X_\tau:\bar\Omega\to M$ where
\begin{equation}
\label{perturbation}
X_\tau (x) = \Xi(\tau, \vartheta(u(x),x)), \quad x\in \bar\Omega. 
\end{equation}
It follows from the Implicit Function Theorem that there exists $\Omega_\tau \subset P$ and $u_\tau:\bar\Omega_\tau\to \mathbb{R}$ such that $\Sigma_\tau$ is the graph of $u_\tau$. Moreover, 
(\ref{Zboundary}) implies that the $\Omega_\tau\subset\Omega$. 
 
Hence given a point $y\in \Sigma$, denote $y_\tau = \Xi(\tau, y)\in \Sigma_\tau$. It follows that there exists $x_\tau\in \Omega_\tau$ such that $y_\tau= \vartheta(u_\tau(x_\tau), x_\tau)$. Then we denote 
by $\hat y_\tau = \vartheta(u(x_\tau), x_\tau)$ the point in $\Sigma$ in the flow line of $Y$ passing through $y_\tau$.  The vertical separation between $y_\tau$ and $\hat y_\tau$ is by definition the function $s(y,\tau)=u_\tau(x_\tau)- u(x_\tau)$.

\begin{lemma}\label{lema1} For any $\tau\in (-\varepsilon, \varepsilon)$, let $A_\tau$ and $H_\tau$ be, respectively,  the Weingarten map and the mean curvature of the hypersurface $\Sigma_\tau$ calculated with respect to the unit normal vector field $N_\tau$ along $\Sigma_\tau$ which satisfies $\langle N_\tau, Y\rangle >0$.  Denote $H=H_0$ and $A=A_0$. If
$\zeta\in C^\infty(\mathcal{O})$ and $T\in \Gamma(T\mathcal{O})$ are defined by
\begin{equation}
Z = \zeta N_\tau + T
\end{equation}
with $\langle T, N_\tau\rangle=0$ then 
\begin{itemize}
\item[i.]
$\frac{\partial s}{\partial\tau}\big|_{\tau=0} = \langle Z, N\rangle  W.$
\item[ii.]
$\bar{\nabla}_Z N\big|_{\tau=0} = -AT-\nabla^{\Sigma}\zeta$
\item[iii.]
$\frac{\partial H}{\partial\tau}\big|_{\tau=0}=\Delta_\Sigma\zeta+(|A|^2+{\rm Ric}_M(N,N))\zeta+\langle\bar\nabla \Psi, Z\rangle,$
\end{itemize}
where $W=\langle Y, N_\tau\rangle^{-1}=(\gamma+|\nabla u_\tau|^2)^{-1/2}$.  The operators $\nabla^\Sigma$ and $\Delta_\Sigma$ are, respectively, the intrinsic gradient operator and the Laplace-Beltrami operator in $\Sigma$ with respect to the induced metric.  Moreover, $\bar\nabla$ and ${\rm Ric}_M$ denote, respectively, the Riemannian covariant derivative and the Ricci tensor in $M$.
\end{lemma}

\noindent \textit{Proof.} (i) Let $(x^i)_{i=1}^n$ a set of local coordinates in $\Omega\subset P$. Differentiating (\ref{perturbation}) with respect to $\tau$ we obtain 
\begin{eqnarray*}
X_{\tau*}\frac{\partial}{\partial\tau} = Z|_{X_\tau} = \zeta N_\tau + T
\end{eqnarray*}
On the other hand differentiating both sides of 
\[
X_\tau(x) =\vartheta(u_\tau(x_\tau), x_\tau)
\]
with respect to $\tau$ we have
\begin{eqnarray*}
X_{\tau*}\frac{\partial}{\partial\tau} &=&\Big( \frac{\partial u_\tau}{\partial \tau}+\frac{\partial u_\tau}{\partial x^i}\frac{\partial x_\tau^i}{\partial \tau}\Big)\vartheta_* Y +\frac{\partial x_\tau^i}{\partial \tau}
\vartheta_* \frac{\partial}{\partial x^i}\\
& = & \frac{\partial u_\tau}{\partial \tau}\vartheta_* Y+\frac{\partial x_\tau^i}{\partial \tau}\Big(\vartheta_* \frac{\partial}{\partial x^i}+\frac{\partial u_\tau}{\partial x^i}\vartheta_* Y\Big)
\end{eqnarray*}
Since the term between parenthesis after the second equality is a tangent vector field in $\Sigma_\tau$ we conclude that
\begin{eqnarray*}
\frac{\partial u_\tau}{\partial \tau}\langle Y, N_\tau\rangle = \langle X_{\tau*}\frac{\partial}{\partial\tau}, N_\tau\rangle = \zeta 
\end{eqnarray*}
from what follows that 
\[
\frac{\partial u_\tau}{\partial \tau} = \zeta W
\]
and
\begin{eqnarray}
\frac{\partial s}{\partial\tau} = \frac{\partial }{\partial\tau} (u_{\tau}-u) = \frac{\partial u_{\tau}}{\partial\tau} = \zeta W.\nonumber
\end{eqnarray}

\noindent (ii) Now we have
\begin{eqnarray}
& & \langle\bar{\nabla}_{Z}N_\tau,X_*\partial_i\rangle = -\langle N_\tau,\bar{\nabla}_{Z}X_*\partial_i\rangle= -\langle N_\tau,\bar{\nabla}_{X_*\partial_i} Z\rangle= -\langle N_\tau,\bar{\nabla}_{X_*\partial_i} (\zeta N+T)\rangle\nonumber\\
& & \,\, =  -\langle N_\tau,\bar{\nabla}_{X_*\partial_i} T\rangle-\langle N_\tau,\bar{\nabla}_{X_*\partial_i} \zeta N_\tau\rangle= -\langle A_\tau T, X_*\partial_i\rangle- \langle\nabla^{\Sigma}\zeta, X_*\partial_i\rangle,\nonumber
\end{eqnarray}
for any $1\le i\le n$. It follows that
\[
\bar{\nabla}_Z N = -AT-\nabla^{\Sigma}\zeta.
\]

\noindent  (iii)  This is a well-known formula whose proof may be found at a number of references (see, for instance, \cite{gerhardt-book}).  \hfill $\square$

\vspace{3mm}

For further reference, we point out that the Comparison Principle \cite{GT} when applied to (\ref{capillary})-(\ref{neumann-condition})  may be stated in geometric terms as follows. Fixed $\tau$, let $x\in \bar\Omega'$ be a point of maximal vertical separation $s(\cdot, \tau)$. If $x$ is an interior point we have
\[
\nabla u_\tau (x,\tau) -\nabla u(x) = \nabla s (x,\tau) = 0,
\]
what implies that the graphs of the functions $u_\tau$ and $u+s(x,\tau)$ are tangent at their common point $y_\tau =\vartheta(u_\tau(x), x)$. Since the graph of $u+s(x, \tau)$ is obtained from $\Sigma$ only by a translation along the flow lines of $Y$ we conclude that the mean curvature of these two graphs are the same at corresponding points. Since the graph of $u+s(x,\tau)$ is locally above the graph of $u_\tau$ we conclude that
\begin{equation}
\label{comparison-int}
H(\hat y_\tau)\ge H_\tau (y_\tau). 
\end{equation}
If $x\in \partial\Omega\subset \partial\Omega'$ we have
\[
\langle \nabla u_\tau, \nu\rangle|_{x} - \langle \nabla u, \nu\rangle|_x =  \langle \nabla s, \nu\rangle \le 0
\]
since $\nu$ points toward $\Omega$. This implies that 
\begin{equation}
\label{comparison-bdry}
\langle N, \nu\rangle|_{y_\tau} \ge \langle N, \nu\rangle|_{\hat y_\tau}  
\end{equation}

\subsection{Interior gradient estimate}
\label{section-int}

\begin{proposition}\label{interior}
Let $B_R(x_0)\subset \Omega$ where $R<{\rm inj}P$. Then there exists a constant $C>0$ depending on $\beta, C_\Psi, \Omega$ and $K$  such that
\begin{equation}
|\nabla u(x)|\le C\frac{R^2}{R^2 -d^2(x)},  
\end{equation}
where $d={\rm dist}(x_0, x)$ in $P$. 
\end{proposition}

\noindent \emph{Proof.} Fix $\Omega'= B_R(x_0)\subset \Omega$. We consider the vector field $Z$ given by 
\begin{equation}
\label{Zint}
Z=\zeta N, 
\end{equation}
where $\zeta$ is a function to be defined later. Fixed $\tau\in [0, \varepsilon)$,  let $x\in B_R(x_0)$ be a point where the vertical separation $s(\cdot, \tau)$ attains a maximum value. 

If $y=\vartheta(u(x), x)$ it follows that  
\begin{equation}
H_\tau (y_\tau) - H_0(y) = \frac{\partial H_\tau}{\partial\tau}\Big|_{\tau=0}\tau + o(\tau). 
\end{equation}
However the Comparison Principle implies that $H_0(\hat y_\tau)\ge H_\tau (y_\tau)$. Using Lemma \ref{lema1} ($iii$) we conclude that
\begin{eqnarray*}
H_0(\hat y_\tau)- H_0(y) \ge  \frac{\partial H_\tau}{\partial\tau}\Big|_{\tau=0}\tau + o(\tau)= (\Delta_\Sigma\zeta+ |A|^2\zeta + \textrm{Ric}_M(N,N)\zeta)\tau + o(\tau). 
\end{eqnarray*}
Since $\hat y_\tau = \vartheta (-s(y,\tau), y_\tau)$ we have
\begin{eqnarray}
\label{dd}
\frac{d\hat y_\tau}{d\tau}\Big|_{\tau=0} =-\frac{ds}{d\tau}\vartheta_{*}\frac{\partial}{\partial s}+\frac{\partial y_\tau^i}{\partial\tau}\vartheta_{*}\frac{\partial}{\partial x^{i}}
=- \frac{ds}{d\tau} Y + \frac{d y_\tau}{d\tau}\Big|_{\tau=0}=-\frac{ds}{d\tau} Y + Z(y).
\end{eqnarray}
Hence using  Lemma \ref{lema1} ($i$) and (\ref{Zint}) we have
\begin{equation}
\label{dtau}
\frac{d\hat y_\tau}{d\tau}\Big|_{\tau=0}=-\zeta WY+\zeta N.
\end{equation}
On the other hand for each $\tau\in (-\varepsilon, \varepsilon)$ there exists a smooth $\xi: (-\varepsilon, \varepsilon)\to TM$ such that
\[
\hat y_\tau = \exp_y \xi(\tau).
\]
Hence we have
\begin{eqnarray}
\frac{d\hat y_\tau}{d\tau}\Big|_{\tau=0} =\xi'(0).\nonumber
\end{eqnarray}
With a slight abuse of notation we denote $\Psi(s,x)$ by $\Psi(y)$ where $y=\vartheta(s,x)$.
It results that
\begin{equation*}
H_0(\hat y_\tau)- H_0(y) = \Psi(x_\tau, u(x_\tau)) - \Psi(x, u(x)) = \Psi(\exp_y \xi_\tau)-\Psi(y)= \langle \bar\nabla\Psi|_y, \xi'(0)\rangle \tau + o(\tau).  
\end{equation*}
However 
\begin{eqnarray}
\langle\bar\nabla\Psi, \xi'(0)\rangle =\zeta \langle \bar\nabla\Psi, N-WY \rangle= -\zeta W\frac{\partial\Psi}{\partial s}+\zeta\langle \bar\nabla\Psi, N\rangle.
\end{eqnarray}
We conclude that 
\begin{equation*}
-\zeta W\frac{\partial\Psi}{\partial s}\tau+\zeta\langle \bar\nabla\Psi, N\rangle \tau + o(\tau) \ge   (\Delta_\Sigma\zeta+ |A|^2\zeta + \textrm{Ric}_M(N,N)\zeta)\tau + o(\tau).
\end{equation*}
Suppose that 
\begin{equation}
W(x) > \frac{C+|\bar\nabla\Psi|}{\beta}
\end{equation}
for a constant $C>0$ to be chosen later.  Hence we have 
\begin{equation*}
(\Delta_\Sigma\zeta+ \textrm{Ric}_M(N,N)\zeta)\tau  + C\zeta \tau \le  o(\tau). 
\end{equation*}
Following \cite{calle} and \cite{korevaar} we choose
\[
\zeta = 1-\frac{d^2}{R^2},
\]
where $d=\textrm{dist}(x_0, \cdot)$.  It follows that
\begin{eqnarray*}
\nabla^\Sigma\zeta  = -\frac{2d}{R^2}\nabla^\Sigma d
\end{eqnarray*}
and
\begin{eqnarray*}
\Delta_\Sigma \zeta = -\frac{2d}{R^2}\Delta_\Sigma d -\frac{2}{R^2}|\nabla^\Sigma d|^2
\end{eqnarray*}
However using the fact that $P$ is totally geodesic and that $[Y,\bar\nabla d]=0$ we have
\begin{eqnarray}
& & \Delta_{\Sigma}d=\Delta_M d-\langle\bar\nabla_{N}\bar\nabla d,N\rangle + nH\langle\bar\nabla d,N\rangle\nonumber\\
& &\,\, = \Delta_P d -\langle \nabla_{\frac{\nabla u}{W}} \nabla d, \frac{\nabla u}{W}\rangle -\gamma^2\langle Y, N\rangle^2 \langle \bar\nabla_Y \bar\nabla d, Y\rangle+ nH\langle\bar\nabla d,N\rangle\nonumber  \nonumber
\end{eqnarray}
Let $\pi:M\to P$ the projection defined by $\pi(\vartheta(s,x))=x$. Then
\[
\pi_* N = -\frac{\nabla u}{W}.
\]
We denote
\[
\pi_*N^\perp = \pi_* N -\langle \pi_* N, \nabla d\rangle\nabla d.
\]
If $\mathcal{A}_d$ and $\mathcal{H}_d$ denote, respectively, the Weingarten map and the mean curvature of the geodesic ball $B_d(x_0)$ in $P$ we conclude that
\begin{eqnarray}
& & \Delta_{\Sigma}d= n\mathcal{H}_d -\langle \mathcal{A}_d(\pi_* N^\perp), \pi_*N^\perp\rangle +\gamma\langle Y, N\rangle^2 \kappa+ nH\langle\bar\nabla d,N\rangle. \nonumber
\end{eqnarray}
where
\[
\kappa = -\gamma\langle \bar\nabla_Y \bar\nabla d, Y\rangle
\]
is the principal curvature of the Kiling cylinder over $B_d(x_0)$ relative to the principal direction $Y$. Therefore we have
\[
|\Delta_\Sigma d|\le C_1(C_\Psi, \sup_{B_R(x_0)}(\mathcal{H}_d+\kappa), \sup_{B_R(x_0)}\gamma )
\]
in $B_R(x_0)$. Hence setting
\[
C_2 = \sup_{B_R(x_0)}\textrm{Ric}_M
\]
we fix
\begin{equation}
\label{C}
C =\max\{2(C_1+C_2), \sup_{\mathbb{R}\times\Omega} |\bar\nabla \Psi|\}.
\end{equation}
With this choice we conclude that 
\[
C\zeta \le \frac{o(\tau)}{\tau},
\]
a contradiction.  This implies that 
\begin{equation}
W(x) \le \frac{C-|\bar\nabla \Psi|}{\beta}.
\end{equation}
However 
\[
\zeta(z) W(z) + o(\tau) = s(X(z), \tau) \le s(X(x), \tau) = \zeta(x) W(x)+o(\tau),
\]
for any $z\in B_R(x_0)$. It follows that
\[
W(z) \le \frac{R^2-d^2(z)}{R^2-d^2(x)} W(x) + o(\tau) \le \frac{R^2}{R^2-d^2(x)} \frac{C-|\bar\nabla \Psi|}{\beta}+o(\tau) \le \widetilde C \frac{R^2}{R^2-d^2(x)},
\]
for very small $\varepsilon>0$. This finishes the proof of the proposition.  \hfill $\square$

\begin{remark}
\label{sphere}
If $\Omega$ satisfies the interior sphere condition for a uniform radius $R>0$ we conclude that 
\begin{equation}
W(x)\le \frac{C}{d_\Gamma(x)},
\end{equation}
for $x\in \Omega$, where $d_\Gamma(x) ={\rm dist}(x, \Gamma)$.
\end{remark}

\subsection{Boundary gradient estimates} 

Now we establish boundary gradient estimates using other local perturbation of the graph which this time has also tangential components.

\begin{proposition}\label{boundary} Let $x_0\in P$ and $R>0$ such that $3R <{\rm inj}P$. Denote by $\Omega'$ the subdomain $\Omega \cap B_{2R}(x_0)$. Then there exists a positive constant $C=C(R, \beta, \beta', C_\Psi, C_\Phi, \Omega, K)$ such that
\begin{equation}
W(x) \le C,
\end{equation}
for all $x\in \overline\Omega'$. 
\end{proposition}

\noindent \emph{Proof.} Now we consider the subdomain $\Omega'=\Omega\cap B_{R}(x_0)$.  We define
\begin{equation}
Z = \eta N + X,
\end{equation}
where
\[
\eta = \alpha_0  v + \alpha_1 d_\Gamma
\]
and $\alpha_0$ and $\alpha_1$ are positive constants to be chosen and $d_\Gamma$ is a smooth extension of the distance function $\textrm{dist}(\,\cdot\, , \Gamma)$ to $\Omega'$ with $|\nabla d_\Gamma|\le 1$
and
\[
v =4R^2-d^2,
\]
where $d=\textrm{dist}(x_0, \cdot)$. 
Moreover
\[
X = \alpha_0\Phi (v\nu-d_\Gamma\nabla v). 
\]
In this case we have
\begin{eqnarray*}
\zeta = \eta +\langle X, N\rangle = \alpha_0  v + \alpha_1 d_\Gamma + \alpha_0\Phi (v\langle N, \nu\rangle-d_\Gamma\langle N, \nabla v\rangle).
\end{eqnarray*}
Fixed $\tau\in [0,\varepsilon)$, let $x\in\bar\Omega'$ be a point where the maximal vertical separation between $\Sigma$ and $\Sigma_\tau$ is attained. We first suppose that $x\in \textrm{int}(\partial\Omega'\cap \partial\Omega)$.  In this case denoting $y_\tau =\vartheta (u_\tau(x), x)\in \Sigma_\tau$ and $\hat y_\tau=\vartheta(u(x), x)\in \Sigma$ it follows from the Comparison Principle that
\begin{equation}
\langle N_\tau, \nu\rangle|_{y_\tau}\ge \langle N, \nu\rangle|_{\hat y_\tau}. 
\end{equation}
Notice that $\hat y_\tau \in \partial\Sigma$. Moreover since $Z|_{K\cap\mathcal{O}}$ is tangent to $K$ there exists $y\in \partial\Sigma$ such that 
\[
y = \Xi (-\tau, y_\tau).
\]
We claim that 
\begin{equation}
\label{der-1}
|\langle \bar\nabla \langle N_\tau, \nu\rangle, \frac{dy_\tau}{d\tau}\big|_{\tau=0}\rangle| \le \alpha_1 (1-\Phi^2) +\widetilde C\alpha_0
\end{equation}
for some positive constant $\widetilde C=C(C_\Phi, K, \Omega, R)$.

Hence (\ref{neumann-condition}) implies that
\begin{eqnarray*}
\langle N, \nu\rangle|_{\hat y_\tau} - \langle N, \nu\rangle|_{y} = \Phi(\hat y_\tau) - \Phi(y) = \tau \langle \bar\nabla \Phi, \frac{d\hat y_\tau}{d\tau}\big|_{\tau=0}\rangle+ o(\tau).
\end{eqnarray*}
Therefore 
\begin{eqnarray*}
\langle N, \nu\rangle|_{y_\tau} - \langle N, \nu\rangle|_{y} \ge \tau \langle \bar\nabla \Phi, \frac{d\hat y_\tau}{d\tau}\big|_{\tau=0}\rangle+ o(\tau).
\end{eqnarray*}
On the other hand we have
\begin{eqnarray*}
\langle N, \nu\rangle|_{y_\tau} - \langle N, \nu\rangle|_{y} = \tau \langle \bar\nabla \langle N, \nu\rangle,  \frac{dy_\tau}{d\tau}\big|_{\tau=0}\rangle+ o(\tau).
\end{eqnarray*}
We conclude that 
\begin{equation*}
\label{ineq-fund}
\tau \langle \bar\nabla \langle N, \nu\rangle,  \frac{dy_\tau}{d\tau}\big|_{\tau=0}\rangle \ge \tau \langle \bar\nabla \Phi, \frac{d\hat y_\tau}{d\tau}\big|_{\tau=0}\rangle+ o(\tau).
\end{equation*}
Hence we have
\begin{equation*}
\label{ineq-fund}
\alpha_1 (1-\Phi^2)\tau +\widetilde C\alpha_0\tau \ge  \tau \langle \bar\nabla \Phi, \frac{d\hat y_\tau}{d\tau}\big|_{\tau=0}\rangle+ o(\tau).
\end{equation*}
It follows from (\ref{dd}) that 
\begin{equation*}
\label{ineq-fund2}
\alpha_1 (1-\Phi^2) +\widetilde C\alpha_0 \ge   -\zeta W \langle \bar\nabla \Phi,Y\rangle+ \zeta\langle \bar\nabla \Phi, N\rangle+  o(\tau)/\tau.
\end{equation*}
Since 
\[
\langle \bar\nabla\Phi, Y\rangle =\frac{\partial\Phi}{\partial s}\le 0
\]
we conclude that
\begin{equation}
\label{ineq-fund3}
W(x) \le C(C_\Phi, \beta',  K, \Omega, R).
\end{equation}
We now prove the claim.  For that, observe that Lemma \ref{lema1} ($ii$) implies that 
\begin{eqnarray*}
& &\langle N, \nu\rangle|_{y_\tau} - \langle N, \nu\rangle|_{y} = \tau \frac{\partial}{\partial\tau}\Big|_{\tau=0}\langle N_\tau, \nu\rangle|_{y_\tau} + o(\tau) \\
& & \,\, =  \tau (\langle N, \bar\nabla_Z \nu\rangle|_y-\langle AT+\nabla^\Sigma \zeta, \nu\rangle|_y)+o(\tau).
\end{eqnarray*}
Since $Z|_y\in T_y K$ it follows that 
\begin{eqnarray*}
\langle N, \nu\rangle|_{y_\tau} - \langle N, \nu\rangle|_{y} = -\tau (\langle A_K Z,  N\rangle|_y+\langle AT+\nabla^\Sigma \zeta, \nu\rangle|_y)+o(\tau),
\end{eqnarray*}
where $A_K$ is the Weingarten map of $K$ with respect to $\nu$. We conclude that
\begin{equation}
\label{ineq222-2}
-\tau (\langle A_K Z,  N\rangle|_y+\langle AT+\nabla^\Sigma \zeta, \nu\rangle|_y) \ge \tau \langle \bar\nabla \Phi, \frac{d\hat y_\tau}{d\tau}\big|_{\tau=0}\rangle+ o(\tau)
\end{equation}
where 
\[
\nu^T = \nu-\langle N, \nu\rangle N.
\]
We have
\begin{eqnarray*}
\langle \nabla^{\Sigma}\zeta+AT, \nu^T\rangle   = \alpha_0  \langle \nabla v, \nu^T\rangle+\alpha_1\langle\nabla^{\Sigma}d_\Gamma, \nu^T\rangle+\langle\nabla^{\Sigma}\langle X,N\rangle, \nu^T\rangle + \langle AT, \nu^T\rangle.
\end{eqnarray*}
We compute
\begin{eqnarray*}
& & \langle \nabla^{\Sigma}\langle X,N\rangle,\nu^T\rangle =\alpha_0 (v\langle N, \nu\rangle -d_\Gamma\langle N, \nabla v\rangle) \langle \bar\nabla\Phi, \nu^T\rangle \\
& & \,\,\,\, + \alpha_0 \Phi \big(\langle \nabla v, \nu^T\rangle \langle N, \nu\rangle + v (\langle\bar\nabla_{\nu^T} N, \nu\rangle+\langle N, \bar\nabla_{\nu^T}\nu\rangle) - \langle \nabla d_\Gamma, \nu^T\rangle \langle N, \nabla v\rangle\\
& & \,\,\,\, - d_\Gamma
(\langle \bar\nabla_{\nu^T} N, \nabla v\rangle+\langle N, \bar\nabla_{\nu^T}\nabla v\rangle)\big).
\end{eqnarray*}
Hence we have at $y$ that 
\begin{eqnarray*}
& & \langle \nabla^{\Sigma}\langle X,N\rangle,\nu^T\rangle =\alpha_0 (v\Phi -d_\Gamma\langle N, \nabla v\rangle) \langle \bar\nabla\Phi, \nu^T\rangle \\
& & \,\,\,\, + \alpha_0 \Phi \big(\langle \nabla v, \nu^T\rangle \Phi + v (-\langle A\nu^T, \nu^T\rangle+\langle N, \bar\nabla_{\nu}\nu\rangle-\langle N,\nu\rangle\langle N, \bar\nabla_{N}\nu\rangle) \\
& & \,\,\,\,-\langle \nu, \nu^T\rangle \langle N, \nabla v\rangle - d_\Gamma
(-\langle A\nu^T, \nabla v\rangle+\langle N, \bar\nabla_{\nu}\nabla v\rangle-\langle N,\nu\rangle\langle N, \bar\nabla_{N}\nabla v\rangle)\big).
\end{eqnarray*}
Therefore we have
\begin{eqnarray*}
& & \langle \nabla^{\Sigma}\langle X,N\rangle,\nu^T\rangle =\alpha_0 (v\Phi -d_\Gamma\langle N, \nabla v\rangle) \langle \bar\nabla\Phi, \nu^T\rangle \\
& & \,\,\,\, +\alpha_0 \Phi \big(\langle \nabla v, \nu^T\rangle \Phi - v (\langle A\nu^T, \nu^T\rangle+\langle N,\nu\rangle\langle N, \bar\nabla_{N}\nu\rangle) \\
& & \,\,\,\,-\langle \nu, \nu^T\rangle \langle N, \nabla v\rangle + d_\Gamma
(\langle A\nu^T, \nabla v\rangle-\langle N, \bar\nabla_{\nu}\nabla v\rangle+\langle N,\nu\rangle\langle N, \bar\nabla_{N}\nabla v\rangle)\big).
\end{eqnarray*}
It follows that 
\begin{eqnarray*}
& & \langle \nabla^{\Sigma}\zeta + AT, \nu^T\rangle =\langle AT, \nu^T\rangle+\alpha_0  \langle \nabla v, \nu^T\rangle+\alpha_1\langle\nu, \nu^T\rangle \\
& & \,\,\,\, +\alpha_0 (v\Phi -d_\Gamma\langle N, \nabla v\rangle) \langle \bar\nabla\Phi, \nu^T\rangle \\
& & \,\,\,\, +\alpha_0 \Phi \big(\langle \nabla v, \nu^T\rangle \Phi - v (\langle A\nu^T, \nu^T\rangle+\langle N,\nu\rangle\langle N, \bar\nabla_{N}\nu\rangle) \\
& & \,\,\,\,-\langle \nu, \nu^T\rangle \langle N, \nabla v\rangle + d_\Gamma
(\langle A\nu^T, \nabla v\rangle-\langle N, \bar\nabla_{\nu}\nabla v\rangle+\langle N,\nu\rangle\langle N, \bar\nabla_{N}\nabla v\rangle)\big).
\end{eqnarray*}
However
\[
\langle AT, \nu^T\rangle= \langle A\nu^T, X\rangle =\alpha_0 \Phi v\langle A\nu^T, \nu^T\rangle -\alpha_0 \Phi d_\Gamma\langle A\nu^T, \nabla v\rangle.
\]
Hence we have
\begin{eqnarray*}
& & \langle \nabla^{\Sigma}\zeta + AT, \nu^T\rangle =\alpha_0  \langle \nabla v, \nu^T\rangle+\alpha_1\langle\nu, \nu^T\rangle  +\alpha_0 (v\Phi -d_\Gamma\langle N, \nabla v\rangle) \langle \bar\nabla\Phi, \nu^T\rangle \\
& & \,\,\,\, +\alpha_0 \Phi \big(\langle \nabla v, \nu^T\rangle \Phi - v \Phi\langle N, \bar\nabla_{N}\nu\rangle -\langle \nu, \nu^T\rangle \langle N, \nabla v\rangle\\
& & \,\,\,\, - d_\Gamma
(\langle N, \bar\nabla_{\nu}\nabla v\rangle-\langle N,\nu\rangle\langle N, \bar\nabla_{N}\nabla v\rangle)\big).
\end{eqnarray*}
Since $d_\Gamma(y)=0$ we have
\begin{eqnarray*}
& & \langle \nabla^{\Sigma}\zeta + AT, \nu^T\rangle =\alpha_0  \langle \nabla v, \nu^T\rangle+\alpha_1\langle\nu, \nu^T\rangle  +\alpha_0 v\Phi \langle \bar\nabla\Phi, \nu^T\rangle \\
& & \,\,\,\, +\alpha_0 \Phi \big(\langle \nabla v, \nu^T\rangle \Phi - v \Phi\langle N, \bar\nabla_{N}\nu\rangle -\langle \nu, \nu^T\rangle \langle N, \nabla v\rangle\big).
\end{eqnarray*}
Rearranging terms we obtain
\begin{eqnarray*}
& & \langle \nabla^{\Sigma}\zeta + AT, \nu^T\rangle =\alpha_1(1-\langle N, \nu\rangle^2)+\alpha_0  \langle \nabla v, \nu^T\rangle (1+\Phi^2) +\alpha_0 v\Phi \langle \bar\nabla\Phi, \nu^T\rangle \\
& & \,\,\,\, -\alpha_0 \Phi \big(v \Phi\langle N, \bar\nabla_{N}\nu\rangle +(1-\langle N, \nu\rangle^2) \langle N, \nabla v\rangle\big).
\end{eqnarray*}
Therefore there exists a constant $C=C(\Phi, K, \Omega, R)$ such that 
\begin{equation}
\label{est-1}
|\langle \nabla^{\Sigma}\zeta + AT, \nu^T\rangle|\le \alpha_1 (1-\Phi^2) +C\alpha_0.
\end{equation}
Since $d_\Gamma(y)=0$ it holds that
\[
|\langle A_K Z, N\rangle| = |A_K| |Z|\le |A_K|(\eta +|X|) \le 4R^2\alpha_0|A_K|(1+\Phi).
\]
from what we conclude that 
\begin{equation}
\label{der-1}
|\langle \bar\nabla\langle N_\tau, \nu\rangle,  \frac{dy_\tau}{d\tau}\big|_{\tau=0}\rangle| \le \alpha_1 (1-\Phi^2) +\widetilde C\alpha_0
\end{equation}
for some constant $\widetilde C(C_\Phi, K, \Omega, R)>0$.

Now we suppose that $x\in \overline{\partial\Omega'\cap \Omega}$. In this case, we have $v(x)=0$. Then $\eta=\alpha_1 d_\Gamma$ and 
\[
X=-\alpha_0 \Phi d_\Gamma\nabla v
\]
at $x$.  Thus 
\[
\zeta = \eta + \langle X,N\rangle = \alpha_1 d_\Gamma +2\alpha_0 \Phi dd_\Gamma \langle \nabla d, N\rangle.
\]
Moreover we have
\[
W(x) \le \frac{C}{d_\Gamma(x)}
\]
(see Remark \ref{sphere}). It follows that
\begin{eqnarray}
\zeta W \le C(\alpha_1 +2\alpha_0 \Phi d\langle \nabla d, N\rangle)\le C(\alpha_1 +4R\alpha_0 \Phi). 
\end{eqnarray}
We conclude that 
\begin{eqnarray}
W(x) \le C(C_\Phi,  K, \Omega, R).
\end{eqnarray}
Now we consider the case when $x\in \Omega \cap \Omega'$. In this case we have
\begin{eqnarray*}
& &\Delta_\Sigma\zeta =\alpha_0  \Delta_\Sigma v + \alpha_1 \Delta_\Sigma d_\Gamma + \alpha_0\Delta_\Sigma\Phi (v\langle N, \nu\rangle-d_\Gamma\langle N, \nabla v\rangle)\\
& & \,\, +\alpha_0\Phi (\Delta_\Sigma v\langle N, \nu\rangle + v\Delta_\Sigma\langle N, \nu\rangle+2\langle \nabla^\Sigma v, \nabla^\Sigma\langle N, \nu\rangle\rangle-\Delta_\Sigma d_\Gamma\langle N, \nabla v\rangle-d_\Gamma\Delta_\Sigma\langle N, \nabla v\rangle\\
& &\,\,\,\, -2\langle \nabla^\Sigma d_\Gamma, \nabla^\Sigma \langle N, \nabla v\rangle)\\
& & \,\,\,\, +2\alpha_0 \langle \nabla^\Sigma \Phi, \nabla^\Sigma v \langle N, \nu\rangle+v \nabla^\Sigma\langle N, \nu\rangle-\nabla^\Sigma d_\Gamma\langle N, \nabla v\rangle-d_\Gamma\nabla^\Sigma\langle N, \nabla v\rangle\rangle 
\end{eqnarray*}
Notice that given an arbitrary vector field $U$ along $\Sigma$ we have
\begin{equation*}
\langle \nabla^\Sigma \langle N, U\rangle, V\rangle = -\langle AU^T, V\rangle +\langle N, \bar\nabla_V U\rangle,  
\end{equation*}
for any $V\in \Gamma(T\Sigma)$. Here, $U^T$ denotes the tangential component of $U$. 
Hence using Codazzi's equation we obtain 
\begin{eqnarray*}
\Delta_\Sigma\langle N, U\rangle \le  \langle \bar\nabla (nH), U^T\rangle +\textrm{Ric}_M (U^T, N) + C|A|
\end{eqnarray*}
for a constant $C$ depending on $\bar\nabla U$ and $\bar\nabla^2 U$. Hence using (\ref{capillary}) we conclude that 
\begin{eqnarray}
\Delta_\Sigma\langle N, U\rangle \le  \langle \bar\nabla \Psi, U^T\rangle +\widetilde C|A|
\end{eqnarray}
where $\widetilde C$ is a positive constant depending on $\bar\nabla U, \bar\nabla^2 U$ and $\textrm{Ric}_M$. 

We also have
\begin{eqnarray*}
\Delta_\Sigma d_\Gamma & = & \Delta_P d_\Gamma +\gamma \langle\bar\nabla_Y\bar\nabla d, Y\rangle -\langle \bar\nabla_N \bar\nabla d_\Gamma, N \rangle +nH\langle \bar\nabla d_\Gamma, N\rangle  \\
& \le  &  C_0\Psi + C_1,  
\end{eqnarray*}
where $C_0 $ and $C_1$ are positive constants depending on the second fundamental form of the Killing cylinders over the equidistant sets $d_\Gamma = \delta$ for small values of $\delta$.  Similar estimates also hold for $\Delta_\Sigma d$ and then for $\Delta_\Sigma v$.

We conclude that 
\begin{equation}
\Delta_\Sigma \zeta \ge -\widetilde C_0 - \widetilde C_1 |A|, 
\end{equation}
where $\widetilde C_0$ and $\widetilde C_1$ are positive constants depending on $\Omega$, $K$, $\textrm{Ric}_M$, $|\Phi|_2$. 

Now proceeding similarly as in the proof of Proposition \ref{interior}, we observe that Lemma \ref{lema1} ($iii$)  and the Comparison Principle yield
\begin{eqnarray*}
H_0(\hat y_\tau)- H_0(y) \ge  \frac{\partial H_\tau}{\partial\tau}\Big|_{\tau=0}\tau + o(\tau)= (\Delta_\Sigma\zeta+ |A|^2\zeta + \textrm{Ric}_M(N,N)\zeta)\tau +\tau\langle \bar\nabla\Psi, T\rangle+ o(\tau). 
\end{eqnarray*}
However
\begin{equation*}
H_0(\hat y_\tau)- H_0(y) =  \langle \bar\nabla\Psi|_y, \xi'(0)\rangle \tau + o(\tau).  
\end{equation*}
Using (\ref{dd}) we have
\begin{eqnarray*}
\langle\bar\nabla\Psi, \xi'(0)\rangle =\langle \bar\nabla\Psi, Z-\zeta WY \rangle= \langle \bar\nabla\Psi, Z\rangle-\zeta W\frac{\partial\Psi}{\partial s}.
\end{eqnarray*}
We conclude that 
\begin{equation*}
-\zeta W\frac{\partial\Psi}{\partial s}\tau+\zeta\langle \bar\nabla\Psi, N\rangle \tau + o(\tau) \ge   (\Delta_\Sigma\zeta+ |A|^2\zeta + \textrm{Ric}_M(N,N)\zeta)\tau + o(\tau).
\end{equation*}
Suppose that 
\begin{equation}
W > \frac{C+|\bar\nabla\Psi|}{\beta}
\end{equation}
for a constant $C>0$ as in (\ref{C}).  Hence we have 
\begin{equation*}
(\Delta_\Sigma\zeta+|A|^2\zeta+ \textrm{Ric}_M(N,N)\zeta)\tau  + C\zeta \tau \le  o(\tau) 
\end{equation*}
We conclude that 
\begin{equation*}
- C_0 - C_1 |A| + C_2 |A|^2 +C\le \frac{o(\tau)}{\tau},
\end{equation*}
a contradiction. It follows from this contradiction that
\begin{equation}
W(x) \le \frac{C+|\bar\nabla\Psi|}{\beta}.
\end{equation}
Now, proceeding as in the end of the proof of Proposition \ref{interior}, we use the estimate for $W(x)$ in each one of the three cases for obtaining a estimate for $W$ in $\Omega'$.  This finishes the proof of the Proposition.   $\hfill\square$

\section{Proof of  the Theorem \ref{main}}
\label{section-proof}

We use the classical Continuity Method for proving Theorem \ref{main}. For details, we refer the reader to \cite{gerhardt} and \cite{uraltseva-book}. For any $\tau\in [0,1]$ we consider the Neumann boundary problem $\mathcal{N}_\tau$ of finding $u\in C^{3,\alpha}(\bar\Omega)$ such that
\begin{eqnarray}
& & \mathcal{F}[\tau, x,u,\nabla u, \nabla^2 u] = 0,\\
& & \langle \frac{\nabla u}{W}, \nu\rangle + \tau \Phi=0,
\end{eqnarray}
where $\mathcal{F}$ is the quasilinear elliptic operator defined by 
\begin{eqnarray}
\mathcal{F}[x,u,\nabla u, \nabla^2 u]= \textrm{div}\bigg(\frac{\nabla u}{W}\bigg) - \langle \frac{\nabla \gamma}{2\gamma}, \frac{\nabla u}{W}\rangle -\tau\Psi.
\end{eqnarray}
Since the coefficients of the first and second order terms do not depend on $u$ it follows that
\begin{equation}
\label{implicit}
\frac{\partial\mathcal{F}}{\partial u}= -\tau\frac{\partial\Psi}{\partial u} \le -\tau\beta <0.
\end{equation}
We define $\mathcal{I}\subset [0,1]$ as the subset of values of $\tau\in [0,1]$ for which the Neumann boundary problem $\mathcal{N}_\tau$ has a solution. Since $u=0$ is a solution for $\mathcal{N}_0$, it follows that $\mathcal{I}\neq \emptyset$. Moroever, the Implicit Function Theorem (see \cite{GT}, Chapter 17) implies that $\mathcal{I}$ is open in view of (\ref{implicit}). Finally, the height and gradient \emph{a priori} estimates we obtained in Sections \ref{section-height} and \ref{section-gradient} are independent of $\tau\in [0,1]$. This implies that (\ref{capillary}) is uniformly elliptic. Moreover, we may assure the existence of some $\alpha_0 \in  (0,1)$ for which there there exists a constant $C>0$ independent of $\tau$ such that 
\[
|u_\tau|_{1,\alpha_0,\bar\Omega}\le C. 
\]
Redefine $\alpha = \alpha_0$. Thus, combining this fact, Schauder elliptic estimates and the compactness of $C^{3,\alpha_0}(\bar\Omega)$ into $C^3(\bar\Omega)$ imply that $\mathcal{I}$ is closed. It follows that $\mathcal{I}=[0,1]$.

The uniqueness follows from the Comparison Principle for elliptic PDEs. We point out that a more general uniqueness statement  - comparing a nonparametric solution with a general hypersurface with the same mean curvature and contact angle at corresponding points - is also valid. It is a consequence of a flux formula coming from the existence of a Killing vector field in $M$. We refer the reader to \cite{DHL} for further details. 

This finishes the proof of the Theorem \ref{main}.

\vspace{3cm}
\noindent 
Jorge H. Lira\\
Gabriela A. Wanderley\\
Departamento de Matem\'atica \\ Universidade Federal do Cear\'a\\ 
Campus do Pici, Bloco 914\\ Fortaleza, Cear\'a\\ Brazil\\ 60455-760

\end{document}